\documentclass[12pt]{amsart}
\usepackage{amsmath, amssymb, amscd}
\usepackage{hyperref, graphicx}
\usepackage{amsrefs}

\newcommand{\PP}{\ensuremath{\mathbb{P}}}
\newcommand{\ZZ}{\ensuremath{\mathbb{Z}}}

\newcommand{\mO}{\ensuremath{\mathcal{O}}}

\newcommand{\mF}{\ensuremath{\mathcal{F}}}

\newcommand{\mI}{\ensuremath{\mathcal{I}}}

\newcommand{\bx}{\ensuremath{{\bf x}}}
\newcommand{\by}{\ensuremath{{\bf y}}}
\newcommand{\bz}{\ensuremath{{\bf z}}}

\newcommand{\wts}{\widetilde{\Sigma}}

\DeclareMathOperator{\codim}{codim}

\newtheorem{thm}{Theorem}
\newtheorem{cor}[thm]{Corollary}

\newtheorem{prop}[thm]{Proposition}

\theoremstyle{definition}
\newtheorem{conj}[thm]{Conjecture}
 \newtheorem{defin}[thm]{Definition}
\newtheorem{ex}[thm]{Example}

\newtheorem{rmk}[thm]{Remark}

\begin{document}

\title{Equations defining secant varieties: geometry and computation}
\author{Jessica Sidman and Peter Vermeire}
\maketitle
\begin{abstract}
In the 1980's, work of Green and Lazarsfeld \cite{GL1, GL2} helped to uncover the beautiful interplay between the geometry of the embedding of a curve and the syzygies of its defining equations.  Similar results hold for the first secant variety of a curve, and there is a natural conjectural picture extending to higher secant varieties as well.   We present an introduction to the algebra and geometry used in \cite{sidver} to study syzygies of secant varieties of curves with an emphasis on examples of explicit computations and elementary cases that illustrate the geometric principles at work.
\end{abstract}
\section{Introduction}

Intuition about the behavior of equations definining a secant variety embedded in projective space arises from consideration of both algebra and geometry, and our main goal is to bring together some of these ideas.   In this paper we will be concerned with syzygies of secant varieties of smooth curves.  We will begin with the algebraic point of view with the aim of giving the reader tools for computing secant varieties of curves using Macaulay 2 \cite{mac2}.
We then turn to the geometric point of view, based on work of Aaron Bertram \cite{bertram}, which led to the second author's original conjectures on cubic generation of secant ideals and linear syzygies \cite{vermeiresecreg}.  These conjectures were refined and strengthened in \cite{sidver} using Macaulay 2 \cite{mac2}.  Bertram's setup is also used by Ginensky \cite{ginensky} to study determinantal equations for curves and their secant varieties.  Our hope is to make some of the geometric intuition accessible to readers familiar with \cite{geomSyz} and \cite{hartshorne}, and that the examples we discuss will be of help in reading the existing literature.

We want to study the minimal free resolution of the homogeneous coordinate ring of a secant variety.  If $X \subset \PP^n$ is a variety, by which we mean a reduced, but not necessarily irreducible, scheme, then we define its $k$th \emph{secant variety}, denoted $\Sigma_k,$  to be the Zariski closure of the union of the $k$-planes in $\PP^n$ meeting $X$ in at least $k+1$ points.   We often write $\Sigma$ for $\Sigma_1.$  We are primarily interested in the situation in which the maps $\Gamma(X, \mO_X(k)) \to \Gamma(\PP^n, \mO_{\PP^n}(k))$ are surjective, or equivalently that the homogeneous coordinate ring $S_X$ is \emph{normally generated}, as then $S_X \cong \oplus \Gamma(X, \mO_X(k))$ and geometric techniques can be used to study $S_X.$   We will assume throughout that if $X \subset \PP^n$ is a curve, then it is embedded via a complete linear system.

Two examples of notions that have geometric and algebraic counterparts are \emph{Castelnuovo-Mumford regularity}, or \emph{regularity}, and the Cohen-Macaulay property, both of which can be defined algebraically in terms of minimal free resolutions.  We begin algebraically, and let $M$ be a finitely generated graded module over the standard graded ring $S= k[x_0, \ldots, x_n].$  The module $M$ has a minimal free resolution
\[
0 \to \underset{j}{\oplus} S(-j)^{\beta_{n,j} }\to \cdots \to  \underset{j}{\oplus} S(-j)^{\beta_{1,j} }\to  \underset{j}{\oplus} S(-j)^{\beta_{0,j}} \to M \to 0.\]
The computer algebra package Macaulay 2 \cite{mac2} computes minimal free resolutions and displays the \emph{graded Betti numbers} $\beta_{i,j}$ in a Betti table arranged as below
\[
\begin{array}{c|ccccc}
	& 0 & 1 & 2 & \cdots & j\\
	\hline\\
0	& \beta_{0,0} & \beta_{1,1}&\beta_{2,2} & \cdots & \beta_{j,j}\\
1& \beta_{0,1} & \beta_{1,2}& \beta_{2,3} & \cdots & \beta_{j, 1+j}\\

\vdots\\
i& \beta_{0,i} & \beta_{1,i+1}&\beta_{2,i+2} & \cdots & \beta_{j, i+j}\\

\end{array}
\]

\begin{defin}
The \emph{regularity} of  a finitely generated graded module $M$ is the maximum $d$ such that some $\beta_{j,d+j}$ is nonzero.\end{defin}

\begin{ex}[The graded Betti diagram of a curve, Example 1.4 in \cite{sidver}]\label{ex:g2d9I}
For example, we can compute the graded Betti diagram of a curve of genus 2 embedded in $\PP^7$ using Macaulay 2.
\begin{verbatim}
   	           	    0  1  2  3  4  5 6
              total: 1 19 58 75 44 11 2
                  0: 1 .  .  .  .  . .
                  1: . 19 58 75 44 5 .
                  2: .  .  .  .  . 6 2
\end{verbatim}
As the diagram shows that $\beta_{5,7}$ and $\beta_{6,8}$ are nonzero, we see that the regularity of the homogeneous coordinate ring is 2 and the homogeneous ideal of the curve has regularity 3.
We will discuss this computation in greater depth in Example \ref{ex:g2d9}. \end{ex}

The regularity of the geometric counterpart of a finitely generated graded module, a coherent sheaf on projective space, has a geometric definition which originally appeared on pg. 99 of \cite{mumford}.
\begin{defin}
 The \emph{regularity} of a coherent sheaf $\mF$ on $\PP^n$ is  defined to be the infimum of all $d$ such that $H^i(\PP^n, \mF(d-i))=0$ for all $i >0.$  
 \end{defin}
If $M= \oplus_{j \geq 0} \Gamma(\PP^n, \mF(j))$, then the regularity of $M$ is the maximum of the regularity of $\mF$ and zero.  The reader will find a good discussion in Chapter 4 of \cite{geomSyz}.   It is well-known that if $X \subseteq \PP^n$ is a smooth curve of genus $g$ and degree $d \geq 2g+1$ then the regularity of $\mI_X$ is 2 if $X$ is a rational normal curve and 3 otherwise, and the reader may find a nice exposition in \cite{geomSyz}.    A similar result is true for the first secant variety of a smooth curve.

\begin{rmk}
As this article was going to press, we learned from Adam Ginensky and Mohan Kumar that the image of a normal variety over an algebraically closed field under a proper morphism with reduced, connected fibers may fail to be normal.  This invalidates the proof of the normality of $\Sigma$ in Lemma 3.2 in \cite{vermeireidealreg}.   The arguments in \cite{vermeireidealreg} and the subsequent papers \cite{vermeiresecreg, sidver} go through under the additional hypothesis that $\Sigma$ is normal, which we add below in Theorems \ref{thm:reg} and \ref{thm:sidver}.  
\end{rmk}

\begin{thm}[\cite{vermeiresecreg, sidver}]\label{thm:reg}
Let $X \subseteq \PP^n$ be a smooth curve of genus $g$ and degree $d \geq 2g+3.$   If $\Sigma$ is normal, the regularity of $\mI_{\Sigma}$ is 3 if $X$ is a rational normal curve and 5 otherwise.
\end{thm}

Moreover, it is natural to conjecture:

\begin{conj}[\cite{vermeiresecreg, sidver}]\label{conj:reg}
Let $X \subseteq \PP^n$ be a smooth curve of genus $g$ and degree $d \geq 2g+2k+1.$   The regularity of $\mI_{\Sigma_k}$ is $2k+1$ if $X$ is a rational normal curve and $2k+3$ otherwise.
\end{conj}
This conjecture holds for genus 0 curves, as the $k$th secant variety of a rational normal curve has ideal generated by the maximal minors of a matrix of linear forms, and thus the ideal is resolved by an Eagon-Northcott complex.  The result for genus 1 is proved in \cite{fisher, vbH}.

\begin{defin}
We say that a variety $X \subset \PP^n$ is \emph{arithmetically Cohen-Macaulay} if the depth of the irrelevant maximal ideal of $S = k[x_0, \ldots, x_n]$ on $S_X$ is equal to the Krull dimension of $S_X.$  Via the Auslander-Buchsbaum theorem, this is equivalent to saying that the length of a minimal free resolution of $S_X$ is equal to $\codim X.$
\end{defin}
Using the correspondence between local and global cohomology, one can see that this is the same as requiring
 $H^i(\PP^n, \mI_X(k))=0$ for all $0< i \leq \dim X.$  If $X \subset \PP^n$ is a normally generated smooth curve of degree $d$ and genus $g,$ then it is arithmetically Cohen-Macaulay as normal generation implies $H^1(\PP^n, \mI_X(j))=0$ for $j \geq 1$.  (The cohomology groups vanish automatically for $j \leq 0.$) 
 
 The main result of \cite{sidver} is
 \begin{thm}[\cite{sidver}]\label{thm:sidver}
 If $X \subset \PP^n$ is a smooth curve of genus $g$ and degree $d \geq 2g+3$, and $\Sigma$ is normal, then $\Sigma$ is arithmetically Cohen-Macaulay.
 \end{thm}
 
\begin{rmk}
 As we know that the singular locus of $\Sigma$ is the curve $C$, its normality is equivalent to the arithmetically Cohen-Macaulay condition via Serre's condition.  Indeed, Theorem \ref{thm:sidver} holds if we assume normality of $\Sigma.$   In fact, we know that the ideal of a secant variety of a rational normal curve has a resolution given by an Eagon-Northcott complex, and we also know the graded Betti diagram of the secant varieties of elliptic normal curves via \cite{vbH}, so in these two cases, we do know normality.
\end{rmk}

 We conjecture that $\Sigma_k$ is arithmetically Cohen-Macaulay if $d \geq 2g+2k+1$  and hope that we can use cohomology to limit both the number of rows and the number of columns in the graded Betti diagram of $\mI_{\Sigma_k}$ in general. 
The main difficulty in the cohomological program is that our hypotheses are solely in terms of the positivity of a line bundle on a smooth curve $X,$ and we need to prove vanishings in the cohomology of sheaves on its secant varieties which will necessarily have singularities.  

We begin \S 2 by giving the definition of the ideal of a secant variety and discussing how it may be computed via elimination and prolongation.  We then give several examples of ideals of smooth curves and their secant varieties.  In \S 3 we discuss the geometry of the desingularization of the secant varieties of a curve and how Terracini recursion may be used to study them.   We have not made an attempt to survey the vast literature on secant varieties of higher dimensional varieties here, choosing instead to limit our attention to secant varieties of smooth curves.
\bigskip

\noindent {\bf Acknowledgements}
Code for the computation of prolongations used in our examples was written in conjunction with the paper \cite{sidsu}with the help of Mike Stillman.  We are grateful to Seth Sullivant and Mike Stillman for allowing us to include this code here.  The first author is partially supported by NSF grant DMS-0600471 and the Clare Boothe Luce Program and also thanks the organizers of the conference on Hilbert functions and syzygies in commutative algebra held at Cortona in 2007 as well as organizers of the Abel Symposium.  We thank Mohan Kumar and Adam Ginensky for their communications.

\section{Computing secant varieties}
In this section we will describe how secant ideals may be computed via elimination and via prolongation.  In \S 2.1 we will see that the ideal of $\Sigma_k(X)$ can be defined as the intersection of an ideal in $k+1$ sets of variables with a subring corresponding to the original ambient space.  Thus, it is theoretically possible to compute the secant ideal of any variety whose homogeneous ideal can be written down.  However, elimination orders are computationally expensive, so this method will be unwieldy for large examples.  If $X$ is defined by homogeneous forms of the same degree, then the method of prolongation can be used to compute the graded piece of $I(\Sigma_k(X))$ of minimum possible degree.  This computation is fast, and in many cases yields a set of generators of $I(\Sigma_k(X)).$  We will discuss prolongation in \S 2.2 and give a \emph{Macaulay 2} implementation in Appendix A.

\subsection{Secant varieties via elimination}

Let $X \subset \PP^n$ be a variety with homogeneous ideal $I \subset k[\bx]$.  We can define the $k$th secant ideal of $I$ so that it can be computed via elimination.  We work in a ring with $k+1$ sets of indeterminates $\by_i = (y_{i,0}, \ldots, y_{i,n})$ and let $I(\by_i)$ denote the image of the ideal $I$ under the ring isomorphism $x_j \mapsto y_{i,j}.$

We define the ideal of the \emph{ruled join} of $X$ with itself $k$ times as in Remark 1.3.3 in \cite{FOV}:
\[
J = I(\by_1) + \cdots + I(\by_{k+1}) .\] 
Geometrically, we embed $X$ into $k+1$ disjoint copies of $\PP^n$ in a projective space of dimension $(k+1)(n+1)-1.$  If a point is in the variety defined by $J$, we will see a point of $X$ in each set of $y$-variables.  If we project to the linear space $[y_{1,0}+ \cdots +y_{1,n}: \cdots :y_{k+1,0}+ \cdots +y_{k+1,n}]$ then the points in the image are points whose coordinates are sums of $k+1$ points of $X.$  

In practice, we make the change of coordinates which is the identity on the first $k$ sets of variables and is defined by $y_{k+1,j} \mapsto y_{k+1,j}-y_{1,j}- \cdots - y_{k,j}$ on the last set of variables.   This has the effect of sending  $y_{1,j}+\cdots + y_{k+1,j}$ to $y_{k+1,j},$ so that the ideal of the $k$th secant variety is the intersection of $I(\by_1) + \cdots +I(\by_k) +I(\by_{k+1}-\by_1- \cdots - \by_k)$ with $k[y_{k+1,0}, \ldots, y_{k+1,n}].$

\begin{ex}[The secant variety of two points in $\PP^2$.]
Consider $X = \{ [1:0:0], [0:1:0]\}$ with defining ideal $I = \langle x_0x_1, x_2 \rangle.$  Using the definition above, the ideal of the join is
\[
J = \langle y_{1,0}y_{1,1}, y_{1,2} \rangle + \langle y_{2,0}y_{2,1}, y_{2,2} \rangle \]
Under the change of coordinates $y_{2,j} \mapsto y_{2,j} - y_{1,j}$ we have
\[
\tilde{J}= \langle y_{1,0}y_{1,1}, y_{1,2} \rangle + \langle (y_{2,0}-y_{1,0})(y_{2,1}-y_{1,1}), y_{2,2} - y_{1,2} \rangle
\]
The variety $V(\tilde{J})$ consists of points of the form
\[
\begin{split}
[a:0:0:a:b:0], [c:0:0:d:0:0],\\
 [0:e:0:0:f:0], [0:g:0:h:g:0],
 \end{split}\]
where $[a:b], [c:d], [e:f], [g:h] \in \PP^1.$  Eliminating the first three variables projects $V(\tilde{J})$ into $\PP^2,$ and we see that the image of $V(\tilde{J})$ under this projection is the line joining the two points of $X.$

\end{ex}

Sturmfels and Sullivant use a modification of this definition of the secant ideal in which they first define an ideal in the ring $k[\bx, \by_1, \ldots, \by_{k+1}]$, which has $k+2$ sets of variables, and then eliminate.  Using the notation from before they work with secant ideals by first defining $J' = I(\by_1) + \cdots +I(\by_{k+1}) + \langle \by_1+ \cdots + \by_{k+1}-\bx \rangle$ and then computing $I(\Sigma_k(X)) = J' \cap k[\bx].$  Eliminating the $y$-variables produces an ideal in $k[\bx]$ that vanishes on all points in $\PP^n$ that can be written as the sum of $k+1$ points of $X,$ and hence defines the secant ideal of $\Sigma_k(X).$

\begin{ex}[The secant variety of two points in $\PP^2$ revisited.]
Consider $X = \{ [1:0:0], [0:1:0]\}$ with defining ideal $I = \langle x_0x_1, x_2 \rangle.$  Using the definition of Sturmfels and Sullivant, the secant variety of $X$ is
\[
\begin{split}
J' = \langle y_{1,0}y_{1,1}, y_{1,2} \rangle + \langle y_{2,0}y_{2,1}, y_{2,2} \rangle \\
+ \langle y_{1,0}+y_{2,0}-x_0, y_{1,1}+y_{2,1}-x_1, y_{1,2}+y_{2,2}-x_2\rangle.
\end{split}
\]
The variety $V(J')$ consists of points of the form
\[
\begin{split}
[a+b:0:0:a:0:0:b:0:0], [c:d:0:c:0:0:0:d:0],\\ [0:e+f:0:0:e:0:0:f:0],[h:g:0:0:g:0:h:0:0],
\end{split}
\]
where $[a:b], [c:d], [e:f], [g:h] \in \PP^1.$  Eliminating the $y$-variables projects $V(J')$ into $\PP^2,$ and we see that the image of $V(J')$ under this projection is the line joining the two points of $X.$
\end{ex}

From the point of view of computation, using the first definition of a secant ideal is probably better as it involves computing an elimination ideal in an ambient ring with fewer variables.   The advantage of the second definition reveals itself in proofs, especially involving monomial ideals.  

Indeed, as the linear form $y_{1,j}+ \cdots+y_{k+1, j}-x_j$ is in $J',$ we see that $(y_{1,j}+ \cdots +y_{k+1, j})^m$ is equivalent to $x_j^m$ modulo $J'.$  Therefore, a monomial $f(\bx)$ is in $J' \cap k[\bx]$ if and only if $f(\by_1+\cdots+\by_{k+1}) \in J'.$   This is a key observation in Lemma 2.3 of \cite{StS}.

\subsection{Secant varieties via prolongation}

A sufficiently positive embedding of a variety $X$ has an ideal generated by quadrics.  We expect the ideal of $\Sigma_k(X)$ to be generated by forms of degree $k+2.$  There are many ways of seeing that $I(\Sigma_k(X))$ cannot contain any forms of degree less than $k+2.$  This fact was made explicit algebraically by Catalano-Johnson \cite{catalano} who showed that $I(\Sigma_k(X))$ is contained in the $(k+1)$st symbolic power of $I(X)$ and cites an independent proof due to Catalisano.  The stronger statement, that $I(\Sigma_k(X))_{k+2} = I(X)^{(k+1)}$ appears in \cite{LM1, LM2} and a proof of a generalization for ideals generated by forms of degree $d$ is in  \cite{sidsu}.

The connection between symbolic powers and the ideals of secant varieties, at least in the case of smooth curves, is implicit in work of Thaddeus \cite{thaddeus}.  Specifically, in \S 5.3 he constructs a sequence of flips whose exceptional loci are the transforms of secant varieties and then identifies the ample cone at each stage.  The identification of sections of line bundles on the transformed spaces with those on the original in \S 5.2 then provides the connection.  Though not made algebraically explicit, this connection is used in \S 2.12 of \cite{wahl}, is present throughout \cite{vermeireidealreg} and is discussed on pg. 80 of \cite{vermeireSecBir}.

The observation that $I(\Sigma_k(X))$ and the $(k+1)$st symbolic power agree in degree $k+2$ leads to an algorithm for quickly computing all forms of degree $k+2$ in $I(\Sigma_k(X)).$  First, we define the \emph{prolongation} of a vector space $V$ of homogeneous forms of degree $d$ to be the space of a forms of degree $d+1$ whose first partial derivatives are all in $V.$  An easy way to compute the prolongation of $V$ is to compute the vector space $V_i$ of forms of degree $d+1$ formed by integrating the elements of $V$ with respect to $x_i.$  Then $V_1\cap \cdots \cap V_n$ is the prolongation of $V.$  We provide \emph{Macaulay2} code for the computation of prolongations in \emph{Appendix A.}

If $V = I(X)_2,$ then the prolongation of $V$ is $I(\Sigma_1(X))_3.$  The prolongation of $I(\Sigma_k(X))_{k+2}$ is $I(\Sigma_{k+1}(X))_{k+3}.$  As each of these spaces is just the intersection of a set of finite dimensional linear spaces, the vector spaces $I(\Sigma_k(X)_{k+2})$ can be computed effectively in many variables.   

\subsection{Computing the ideal of a smooth curve}

In this section we discuss the examples that we have been able to compute so far, most of which also appear in \cite{sidver}.  Essentially we need a mechanism for writing down the generators of $I(X).$ 
A curve of degree $\geq 4g+3$ has a determinantal presentation by \cite{EKS}, where matrices whose 2-minors generate $I(X)$ for elliptic and hyperelliptic curves are given.   We can also re-embed plane curves into higher dimensional spaces.  We give explicit examples of each class of example below.  We computed the ideals of the secant varieties using the idea of prolongation described in the previous section.  In each case, the projective dimension of the ideal generated by prolongation is equal to the codimension of the secant variety.  Hence, we can deduce that the ideal of the secant variety is generated by the prolongation.

\begin{ex}[Re-embedding a plane curve with nodes: $g=2, d=5$, Example 1.5 in \cite{sidver}]\label{ex:g2d9}
Suppose we have a plane quintic with 4 nodes.  If we blow up the nodes we have a smooth curve of genus 2.  We provide Macaulay 2 \cite{mac2} code below that shows how to compute the ideal of the embedding of such a curve in $\PP^7.$  This method of finding the equations of a smooth curve is due to F. Schreyer and was suggested to the first author by D. Eisenbud.

\begin{verbatim}
--The homogeneous coordinate ring of P^2.
S= ZZ/32003[x_0..x_2]

--These are the ideals of the 4 nodes
I1 = ideal(x_0, x_1)
I2 = ideal(x_0, x_2)
I3 = ideal(x_1, x_2)
I4 = ideal(x_0-x_1, x_1-x_2)

--Forms in I vanishing  twice at each of our chosen nodes
I = intersect(I1^2, I2^2, I3^2, I4^2);  

--The degree 5 piece of I.
M = flatten entries gens truncate(5, I)

--The target of the rational map given by M.
R = ZZ/32003[y_0..y_8]

--The rational map given by M and the ideal of its image.
f = map(S, R, M)
K = ker f

--A random linear change of coordinates on P^8.
g = map(R, R, {random(1, R), random(1, R),random(1, R),
random(1, R),random(1, R),random(1, R),random(1, R),
random(1, R),random(1, R)})

--Add in the element y_8 and then eliminate it.  
--J = ideal of the cone over our plane quintic in P^7.
J =eliminate(y_8, g(K)+ideal(y_8));
\end{verbatim}
The graded Betti diagram of the ideal of the curve is Example \ref{ex:g2d9I} and also Example 1.5 in \cite{sidver}.  Its secant ideal has graded Betti diagram given below.
\begin{verbatim}
                          0  1  2 3 4
                   total: 1 12 16 8 3
                       0: 1  .  . . .          
                       1: .  .  . . .          
                       2: . 12 16 . .          
                       3: .  .  . 4 .          
                       4: .  .  . 4 3
\end{verbatim}

\end{ex}

\begin{ex}[A determinantal curve: $g = 2, d=12$, Example 4.8 in \cite{sidver}]\label{ex:det}
Following \cite{EKS}, we can write down matrices whose $2 \times 2$ minors generate the ideal of a hyperelliptic curve.  We give such a matrix below for a curve with $g=2$ and $d=12.$

\[\begin{pmatrix}{x}_{0}&      {x}_{1}&      {x}_{2}&      {x}_{3}&      {y}_{0}\\      {x}_{1}&      {x}_{2}&      {x}_{3}&      {x}_{4}&      {y}_{1}\\      {x}_{2}&      {x}_{3}&      {x}_{4}&      {x}_{5}&      {y}_{2}\\      {x}_{3}&      {x}_{4}&      {x}_{5}&      {x}_{6}&      {y}_{3}\\      {y}_{0}&      {y}_{1}&      {y}_{2}&      {y}_{3}&      {x}_{1}+{x}_{2}+{x}_{3}+{x}_{4}+{x}_{5}\\      \end{pmatrix}\]
\bigskip

The varieties $\Sigma_k$, $k=0,1,2$, have the graded Betti diagrams below.  Notice that the index of the final row of each diagram is $2k,$ indicating that the homogenous coordinate ring of $\Sigma_k$ has regularity $2k$ and that the ideal of $\Sigma_k$ has regularity $2k+1$ as predicted by Conjecture~\ref{conj:reg}.  Moreover, our curve sits in $\PP^{10},$ and $\dim \Sigma_k = 2k+1,$ so the index of the last column is the codimension of $\Sigma_k,$ indicating that each variety is arithmetically Cohen-Macaulay.  We can also see that $I(\Sigma_k)$ is generated in degree $k+2$ and has linear syzygies up to stage $p$ where $12 = 2g+2k+1+p$ and that $\beta_{9-2k, 11} = \binom{g+k}{k+1}$ as in Conjecture 1.4 in \cite{sidver}.

\begin{verbatim}

                0  1   2   3   4   5   6   7  8 9
         total: 1 43 222 558 840 798 468 147 17 2          
             0: 1  .   .   .   .   .   .   .  . .          
             1: . 43 222 558 840 798 468 147  8 .          
             2: .  .   .   .   .   .   .   .  9 2
             \end{verbatim}

\pagebreak
\begin{verbatim}
               0  1   2   3   4   5  6 7
        total: 1 70 283 483 413 155 14 3         
            0: 1  .   .   .   .   .  . .         
            1: .  .   .   .   .   .  . .         
            2: . 70 283 483 413 155  . .         
            3: .  .   .   .   .   .  7 .         
            4: .  .   .   .   .   .  7 3

              0  1  2  3  4 5
       total: 1 41 94 61 11 4         
           0: 1  .  .  .  . .         
           1: .  .  .  .  . .         
           2: .  .  .  .  . .         
           3: . 41 94 61  . .         
           4: .  .  .  .  . .         
           5: .  .  .  .  6 .         
           6: .  .  .  .  5 4

\end{verbatim}
\end{ex}

\begin{ex}[A Veronese re-embedding of a plane curve: $g=3, d=12$]\label{ex:ver}
Let $X$ be a smooth plane curve of degree 4 and genus 3.  Re-embedding this curve via the degree 3 Veronese map we have a curve of degree 12 in $\PP^9.$  Below we give the graded Betti diagram of the curve and its first two secant varieties.
\begin{verbatim}

             0  1   2   3   4   5  6  7 8
     total: 1 33 144 294 336 210 69 16 3
          0: 1  .   .   .   .   .  .  . .
          1: . 33 144 294 336 210 48  . .
          2: .  .   .   .   .   . 21 16 3

             0  1   2   3  4  5 6
      total: 1 38 108 102 43 18 6
          0: 1  .   .   .  .  . .
          1: .  .   .   .  .  . .
          2: . 38 108 102 10  . .
          3: .  .   .   . 30  . .
          4: .  .   .   .  3 18 6
	  \end{verbatim}
	  \pagebreak
	  \begin{verbatim}
             0 1  2  3  4
      total: 1 8 23 26 10
          0: 1 .  .  .  .
          1: . .  .  .  .
          2: . .  .  .  .
          3: . 8  .  .  .
          4: . .  6  .  .
          5: . . 16 10  .
          6: . .  1 16 10
\end{verbatim}

\end{ex}

Examples \ref{ex:g2d9}, \ref{ex:det} and \ref{ex:ver} together suggest an additional conjecture, that row $2k$ of the Betti diagram of $\Sigma_k$ has precisely $g$ nonzero elements.
\section{Secant varieties as vector bundles: Terracini recursion}
The definition of a secant variety as the Zariski closure of the union of secant lines of $X$ does not lend itself to thinking about a secant variety geometrically.   A more elegant point of view is to realize that a secant line to a projective variety $X$ is just the span of a length two subscheme of $X$.  Thus we should think of $\Sigma(X)$ as the image of a $\PP^1$-bundle over the space of length two subschemes of X, $\operatorname{Hilb}^2X$.  One nice consequence of this point of view is that $\operatorname{Hilb}^2X=\operatorname{Bl}_{\Delta}(X\times X)/S_2$ is smooth as long as $X$ is, and we obtain a geometric model for the secant variety on which we can apply standard cohomological techniques.  Perhaps more importantly, this $\PP^1$-bundle can be constructed explicitly via blowing up.  

Under mild hypotheses on the positivity of the embedding of $X,$ the blowup of $\PP^n$ at $X$ produces a desingularization of $\Sigma(X)$ as a $\PP^1$ bundle over $\operatorname{Hilb}^2X$.   Thinking of this bundle embedded inside the blowup of $\PP^n$ at $X$ we can also examine how it meets the exceptional divisor.   What we see above a point $p \in X$ is a $\PP^{n-2}$ that meets the proper transform of $\Sigma(X)$ in the projection of $X$ into $\PP^{n-2}$ from the tangent space at $p.$

In this section we provide explicit examples illustrating this point of view.  In \S 3.1 we illustrate the geometry of the blowups desingularizing secant varieties when $X$ is a set of 5 points in $\PP^3.$ Although the secant varieties of rational normal curves are well-understood algebraically, we discuss them here as we may make explicit computations and give proofs which highlight the main ideas used in the more general case in \cite{vermeiresecreg} but are much simpler.  We make computations with rational normal curves of degrees 3 and 4 in \S 3.2 and \S 3.3.   In \S 3.4 we discuss how to think about the cohomology along the fibers of these blowups and how cohomology may be used to show that the secant variety is projectively normal.

\subsection{Secant varieties of points}

Let $X\subset \PP^n$ be a finite set of points in linearly general position.  We will analyze the successive blowups of $X$ and the proper transforms of its secant varieties, moving up one dimension at each stage.  It is instructive to consider the geometry of the blowups of a finite set of points and its secant varieties as  we can easily restrict our attention to the picture above a single point.   We think the general picture of the geometry of the blowups will become transparent if we illustrate the construction in a concrete example.  Although this computation may seem quite special, if $X$ consists of $n+2$ points in linearly general position, then a result of Kapranov \cite{kapranov} tells us that the sequence of blowups actually gives a realization of $\overline{M}_{0,n}.$

\subsubsection{Points in $\PP^3$}

Let $B_0  = \PP^3$ and let $\Sigma_0=X$ be a set of 5 linearly general points denoted $p_1, p_2, p_3, p_4, p_5$.  Let $B_1$ be the blowup of $B_0$ at $\Sigma_0.$  We let $E_i$ denote the exceptional divisor above $p_i$ and $\wts_j$ be the proper transform of $\Sigma_j$ for $j = 1,2.$  

\subsubsection{The first recursion}
The exceptional divisor $E_i$ is a $\PP^2$ in which  $\wts_1 \cap E_i$ is a set of 4 points in $\PP^2$ and $\wts_2 \cap E_i$ is the union of lines joining these points in $\wts_1$.  Below we give a diagram depicting the exceptional divisor above a point $p_1$ together with the strict transforms of the span of $p_1$ with two other points $p_2$ and $p_3.$
\[
\includegraphics[width=6cm]{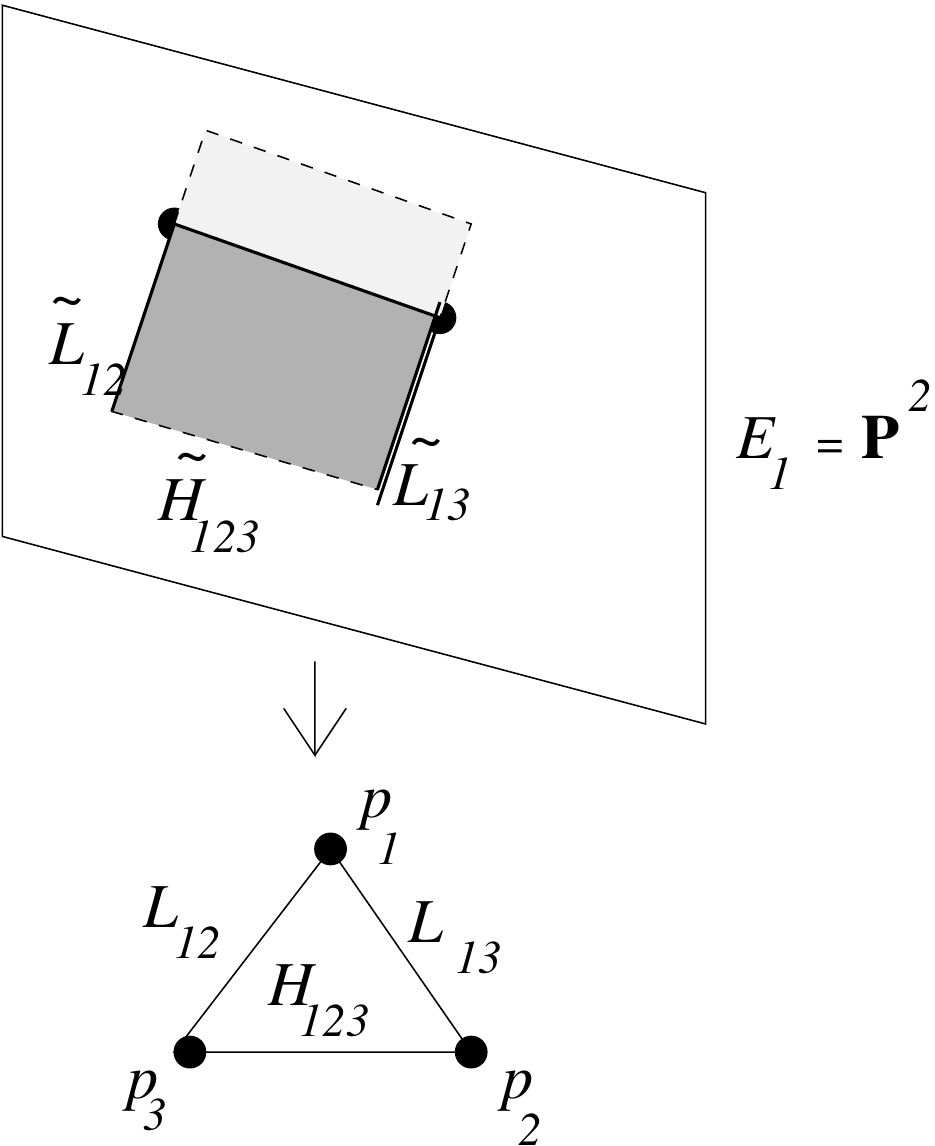}
\]
  The picture in $E_1$ can be found by projecting $\Sigma_1$ and $\Sigma_2$ into $\PP^2$ away from $p_1.$

Taking a more global picture, we see that $\wts_1$ is a smooth variety consisting of the disjoint unions of proper transforms of lines.  However, $\wts_2$ is not smooth as the components of $\wts_2$ intersect the $E_i$ in lines which meet at points.

\subsubsection{The second recursion}
Let us now define $B_2$ to be the blowup of $B_1$ at $\wts_1.$  We will abuse notation and let $E_i$ denote its own proper transform in $B_2.$  
To analyze $B_2$,  it may be helpful to restrict our attention locally to a single point $p = [0:0:0:1]$ and examine the fiber over $p$ after blowing up $p$ and then blowing up the proper transform of a line containing $p.$  Using bihomogeneous coordinates $\bx$ and $\by$ on  $\PP^3 \times \PP^2,$ the blowup of $\PP^3$ is defined by $I(B_p) = \langle x_iy_j-x_jy_i \mid i \neq j \in \{0,1,2\} \rangle.$  Blowing up the proper transform of the line $L$ defined by $\langle x_0, x_1 \rangle$ yields a subvariety of $\PP^3 \times \PP^2 \times \PP^1$ which can be given in tri-homogeneous coordinates $\bx, \by, \bz$ by
\[
I(B_{p,L}) = \langle x_iy_j-x_jy_i \mid i \neq j \in \{0,1,2\} \rangle + \langle x_0z_1-x_1z_0, y_1z_0-y_0z_1\rangle.
\]
To understand the fiber above $p,$ we add the ideal of the point to $I(B_{p,L})$ to get $\langle x_0, x_1, x_2, y_1z_0-y_0z_1 \rangle.$  The $\bx$ coordinates alone cut out $p \times \PP^2 \times \PP^1,$ and the equation in the $\by$ and $\bz$-variables cuts out the blowup of the point $p \times [0:0:1]$ in $p \times \PP^2.$  Thus, we can see that if  blow up a point $p$ and a line containing it, above $p$ we get a $\PP^2$ in which we have blown up one point.  

Turning back to the global picture, we see that when we have blown up $\wts_1,$ the proper transform of $\wts_2$ in $B_2$ is smooth.  Moreover, above each $p_i$ we have a copy of our global picture projected into $\PP^2$ away from a point of $X.$  Indeed, each $E_i$ is a  $\PP^2$ in which we have blown up 4 points.  These 4 points correspond to the intersection of $\wts_1$ with $E_i.$  The intersection of the proper transform of $\Sigma_2$ with $E_i$ consists of the union of the exceptional divisors of the 4 points that are blown up in $E_i.$

\subsection{Blowing up a rational normal curve of degree 3}
We begin with the twisted cubic $X\subset\PP^3$.  In this case the equations of the blowup are easy to write down and we can realize the blowup of $\PP^3$ along $X$ as a $\PP^1$ bundle over $\operatorname{Hilb}^2 \PP^1$ explicitly.  The secant variety $\Sigma(X)=\PP^3$ is smooth, but the embedding is positive enough for $I(X)$ to have linear first syzygies, which is the condition required for the setup in \cite{flip1}.

The three quadrics $x_0x_2-x_1^2, x_0x_3-x_1x_2, x_1x_3-x_2^2,$ which generate $I(X)$ give a rational map $\PP^3 \dashrightarrow \PP^2.$  The blowup of $\PP^3$ at $X$ is the graph of this map  in $\PP^3 \times \PP^2$ with bihomogeneous coordinates $\bx$ and $\by.$   This graph is defined set-theoretically by equations coming from Koszul relations, for example $y_0(x_0x_3-x_1x_2)-y_1(x_0x_2-x_1^2).$  But these relations are generated by the two relations coming from linear syzygies:
\[
x_0y_2-x_1y_1+x_2y_0, x_1y_2-x_2y_1+x_3y_0,
\]
which generate the ideal of the blowup.

Let  $\widetilde{\PP}^3$ denote the blowup of $\PP^3$ along $X$  and $E$ denote the exceptional divisor.   Further, let $q_1$ and $q_2$ denote the restrictions of the projections from $\PP^3 \times \PP^2$ to the first and second factors to $\widetilde{\PP}^3.$  Then $q_1:\widetilde{\PP}^3 \to \PP^3$ is the blowup map and $q_2:\widetilde{\PP}^3 \to \PP^2$ is the morphism induced by $|2H-E|.$  We analyze the fibers of both maps explicitly below.

Above any point $p \in X,$ the fiber of $q_1$ is a $\PP^1$.  For example, if $p = [0:0:0:1],$ then $q_1^{-1}(p)$ is defined by adding $\langle x_0, x_1, x_2 \rangle$ to the ideal of the blowup to get  \[\langle x_0, x_1, x_2, x_3y_0\rangle = \langle x_0, x_1, x_2, x_3 \rangle \cap \langle x_0, x_1, x_2, y_0 \rangle.\]
The first primary component is irrelevant, and the second defines the $\PP^1$ with points $([0:0:0:1],[0:y_1:y_2].)$  

Moreover, we can see from these equations that the blowup of $\PP^3$ along the twisted cubic is a $\PP^1$-bundle over $\PP^2$.  To see $q_2^{-1}([1:0:0]),$ add the ideal $\langle y_1, y_2 \rangle$ to the ideal of the blowup to get the ideal $y_0 \langle x_2, x_3 \rangle.$  This shows that the fiber above $[1:0:0]$ consists of points of the form $([a:b:0:0] , [1:0:0]).$  As a length $n$ subscheme of $\PP^1$ has an ideal generated by a single form of degree $n,$  $\operatorname{Hilb}^n\PP^1=\PP^n$, and so this matches exactly what we expect from the description above.

Choosing a less trivial example, in the next section we will begin to see a recursive geometric picture analogous to what we saw when we blew up a finite set of points and its secant varieties.  Again, we will see the projection of our global picture in the fibers above points of our original variety.  Note that we are projecting to the projectivized normal bundle and hence away from the \emph{tangent} space to a point on our variety.  (In our earlier example, we projected from a single point because a zero-dimensional variety has a zero-dimensional tangent spaces.)

\subsection{Blowing up a rational normal curve of degree 4}
Let $X \subset \PP^4$ be a rational normal curve with defining ideal minimally generated by the six $2\times 2$ minors of the matrix
\[
\begin{pmatrix}
x_0 & x_1 & x_2 & x_3\\
x_1& x_2& x_3& x_4
\end{pmatrix}
\]
These six quadrics are a linear system on $\PP^4$ with base locus $X,$ so they give a rational map $\PP^4 \dashrightarrow \PP^5.$   The closure of the graph of this map in $\PP^4 \times \PP^5$ is $B = B_X(\PP^4),$ $\PP^4$ blown up at $X.$

\subsubsection{The ideal of the blowup}
    Let $R = k[\bx, \by].$  Since $I(B)$ defines a subscheme of $\PP^4 \times \PP^5$ it is a bihomogeneous ideal.  It must contain homogeneous forms in the $y$-variables that define the image of $\PP^4$ in $\PP^5.$  Since each $y_i$ is the image of a quadric vanishing on $X$, the ideal $I(B)$ will also contain bihomogeneous forms that are linear in the $y_i$ corresponding to syzygies.

In our example, we get an ideal with 9 generators by running the Macaulay 2 \cite{mac2} code
\begin{verbatim}
--The coordinate ring of P^4 x P^5.
S = ZZ/32003[x_0..x_4, y_0..y_5]

--The coordinate ring of P^4 with an extra parameter.
--The parameter t ensures that the kernel is bi-homogeneous.
B = ZZ/32003[x_0..x_4,t]

--The map S -> B.
f = map(B,S,{x_0, x_1, x_2, x_3, x_4, 
t*(-x_1^2+x_0*x_2), t*(-x_1*x_2+x_0*x_3),
 t*(-x_2^2+x_1*x_3), t*(-x_1*x_3+x_0*x_4),
 t*(-x_2*x_3+x_1*x_4), t*(-x_3^2+x_2*x_4)})

--The ideal defining the blowup in S
K = ideal mingens ker f
\end{verbatim}
One generator, ${y}_{2}^{2}-{y}_{2} {y}_{3}+{y}_{1}{y}_{4}-{y}_{0} {y}_{5},$ in the $\by$-variables alone cuts out the image of $\PP^4$ in $\PP^5.$  The other 8 generators are constructed from linear syzygies on generators of $I(X)$ as in the previous example.  We can use the ideal to analyze what happens when we take the pre-image of a point $p \in \PP^4$ under the blowup map.  There are two cases depending on whether $p$ is contained in $X.$

Case 1:  Suppose that $p=[0:1:0:0:0] \notin X.$  When we compute $I(p)+I(B)$ and then use local coordinates where $x_1=1,$ we get the ideal $\langle x_0, x_2, x_3, x_4, y_0, y_1, y_2, y_3, y_4\rangle.$  Thus, the pre-image of $p$ is a single point as expected.

Case 2:  Suppose that $p=[1:0:0:0:0] \in X.$  Again, we compute $I(p)+I(B)$ and then use local coordinates where $x_0=1.$  We get the ideal $\langle x_1, x_2, x_3, x_4, y_2, y_4, y_5 \rangle$
which defines a $\PP^2.$

\subsubsection{The intersection of $\wts_1$ with the exceptional divisor}
 To examine what happens to the secant variety of $X$ algebraically, we add the equation of $\Sigma_1(X)$ to $I(B)$.   Above the point $[1:0:0:0:0] $ we have the intersection of the $\PP^2$ with coordinates $[y_0:y_1:y_3]$ with the hypersurface defined by $y_1^2-y_0y_3$.  The conic $y_1^2-y_0y_3$ is the projection of our original curve away from the line with points of the form $[a:b:0:0:0].$ We have $[y_0:y_1:y_3] = [t^2:t^3:t^4],$  which is defined by the given equation.  Schematically, we have
 \[
 \includegraphics[width=4cm]{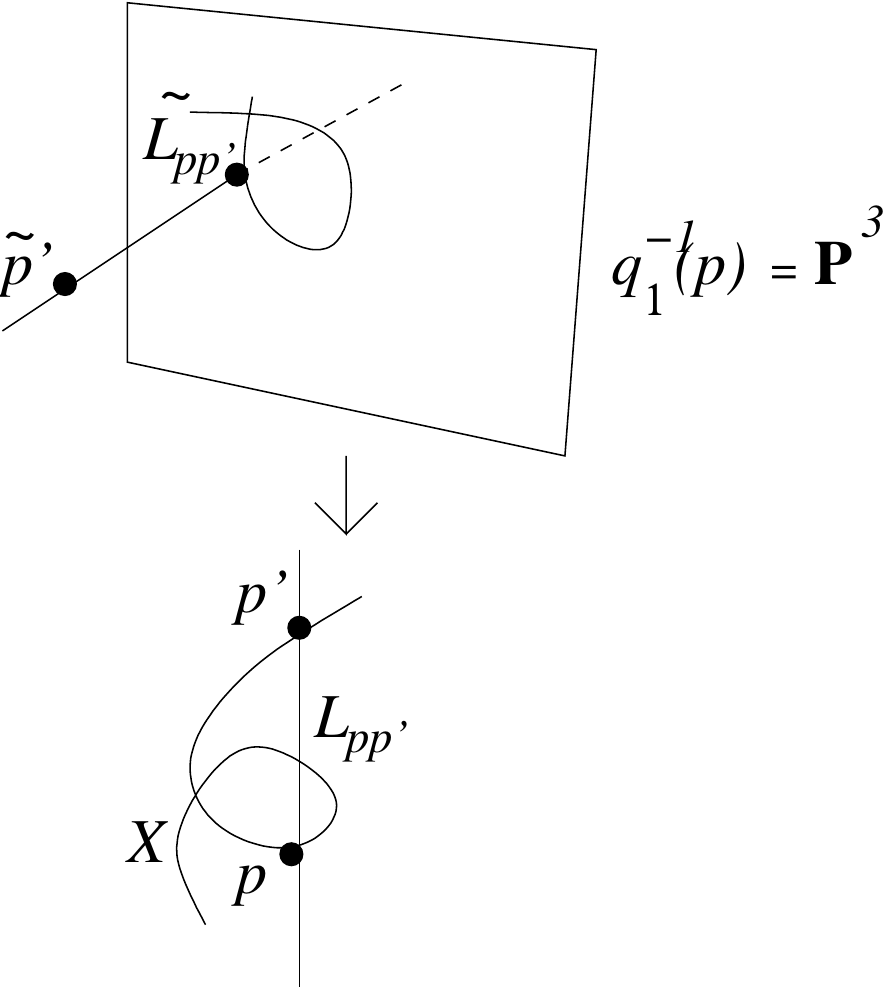}
 \]

\subsection{Cohomology along the fibers}
It was observed by the second author in \cite{vermeireidealreg} that one could use Bertram's \textit{Terracini Recursiveness} to obtain cohomological relationships between different embeddings of the same curve.  For example, let $X$ be a smooth curve embedded in $\PP^n$ by a line bundle $L$ of degree at least $2g+3.$  As discussed above, blowing up $\PP^n$ along $X$ desingularizes $\Sigma_1$, and thus blowing up again along the proper transform of $\Sigma_1$ yields a smooth variety $B_2$
$$B_2 \stackrel{\pi_2}{\rightarrow} B_1  \stackrel{\pi_1}{\rightarrow} \PP^n=\PP\Gamma(X,L).$$
Let $\pi=\pi_1\circ\pi_2$.  

If $x\in X$, then $\pi_1^{-1}(x)\cong\PP^{n-2}=\PP\Gamma(X,L(-2x))$.  By Terracini Recursiveness, $\pi^{-1}(x)$ is the blow up of $\PP\Gamma(X,L(-2x))$ along a copy of $X$ embedded by $L(-2x)$; equivalently $\pi^{-1}(x)$ is precisely what is obtained by the projection $\PP^n\dashrightarrow\PP^{n-2}$ from the line tangent to $X$ at $x$.  In fact, it can be shown that the exceptional divisor of the desingularization $\pi:\widetilde{\Sigma}\rightarrow\Sigma$ is precisely $X\times X$ where the restriction $\pi:X\times X\rightarrow X$ is projection \cite[Lemma 3.7]{flip1}.

Because $B_2$ is obtained from $\PP^n$ by blowing up twice along smooth subvarieties, we know that $\operatorname{Pic}(B_2)=\ZZ H+\ZZ E_1+\ZZ E_2$ where $H$ is the proper transform of a hyperplane section, $E_1$ is the proper transform of the exceptional divisor of the first blow-up $\pi_1$, and $E_2$ is the exceptional divisor of the second blow-up.  Note that we similarly have $\operatorname{Pic}(\pi^{-1}(x))=\ZZ \widetilde{H}+\ZZ \widetilde{E}_1$.  We examine the relationships among line bundles on $\PP^n, B_1,$ and $B_2.$  

Because the generic hyperplane in $\PP^n$ misses $x\in X$, the restriction of $H$ to $\pi^{-1}(x)$ is trivial.  Because $E_1$ is the proper transform of the exceptional divisor of the first blow-up, the restriction of $E_1$ to $\pi^{-1}(x)$ is $-\widetilde{H}$.  Finally, by Terracini Recursiveness, the restriction of $E_2$ to $\pi^{-1}(x)$ is $\widetilde{E}_1$.  Thus a typical effective line bundle on $B_2$ of the form $\mathcal{O}_{B_2}(aH-bE_1-cE_2)$ restricts to $\mathcal{O}_{\pi^{-1}(x)}(b\widetilde{H}-c\widetilde{E}_1)$.

\subsubsection{Regularity and projective normality of the secant variety to a rational normal curve}
In this section will illustrate how to apply the first stage of Bertram's Terracini Recursiveness in the special case where $X$ is a rational normal curve.   Suppose that $X\subset\PP^n$ is a rational normal curve with $L=\mathcal{O}_X(1)=\mathcal{O}_{\PP^1}(n)$ of degree at least $4$ so that $\mathcal{I}_{\Sigma}$ is $3$-regular, and $\Sigma$ is projectively normal. 

It follows \cite[Proposition 9]{vermeiresecreg} from the description of the exceptional divisor of the desingularization $\pi:\widetilde{\Sigma}\rightarrow\Sigma$ above that $\Sigma$ has rational singularities; i.e. $R^i\pi_*\mathcal{O}_{\widetilde{\Sigma}}=0$ for $i>0$.  Thus by Leray-Serre we immediately have $H^i(\widetilde{\Sigma},\mathcal{O}_{\widetilde{\Sigma}}(k))=H^i(\Sigma,\mathcal{O}_{\Sigma}(kH))$ for all $i,k$.

\begin{prop}
Let $X\subset\PP^n$ be a rational normal curve.  Then $\mathcal{I}_{\Sigma}$ is $3$-regular.
\end{prop}

\begin{proof}
We show directly that $H^i(\PP^n,\mathcal{I}_{\Sigma}(3-i))=0$ for $i\geq1$.  As $\Sigma$ is $3$-dimensional, we have only to show the four vanishings $1\leq i\leq4$.  

For each $i \geq 2$ and all $k$ we have $$H^i(\PP^n,\mathcal{I}_{\Sigma}(k))=H^{i-1}(\Sigma,\mathcal{O}_{\Sigma}(k))=H^{i-1}(\widetilde{\Sigma},\mathcal{O}_{\widetilde{\Sigma}}(kH)).$$  

The restriction of $\mathcal{O}_{\widetilde{\Sigma}}(kH))$ to a fiber of the map from $p:\wts \to \PP^2$ is $\mO_{\PP^1}(kH).$  When $i=4$, all cohomology along the fiber vanishes and so in particular we have $H^{3}(\widetilde{\Sigma},\mathcal{O}_{\widetilde{\Sigma}}(-H))=0$.  In general, if $k \geq -1,$ all of the higher cohomology along the fibers vanishes.  This implies that the higher direct image sheaves vanish and that $H^i(\wts, \mO_{\wts}(kH)) = H^i(\PP^2, p_*\mO_{\wts}(kH))).$  When $i=3$, $k =0,$ and  $H^2(\widetilde{\Sigma},\mathcal{O}_{\widetilde{\Sigma}})=H^2(\PP^2, p_*\mO_{\wts}) = H^2(\PP^2,\mathcal{O}_{\PP^2})=0$.

For $i=2$, we blow up $B_1$ along $\widetilde{\Sigma}$ and use Terracini recursiveness.  Consider the sequence
$$0\rightarrow\mathcal{O}_{B_2}(H-E_1-E_2)\rightarrow\mathcal{O}_{B_2}(H-E_2)\rightarrow\mathcal{O}_{E_1}(H-E_2)\rightarrow0.$$
As discussed earlier, the restriction of $\mathcal{O}_{E_1}(H-E_2)$ to a fiber of the flat morphism $E_1\rightarrow X$ is $\mathcal{O}_{\pi^{-1}(x)}(-\widetilde{E}_1)$, but $H^i(\pi_1^{-1}(x),\mathcal{O}_{\pi^{-1}(x)}(-\widetilde{E}_1))=H^i(\PP^{n-2},\mathcal{I}_X)=0$ for $i\geq0$.  Thus it follows that 
\begin{eqnarray*}
H^i(B_2,\mathcal{O}_{B_2}(H-E_1-E_2))&=&H^i(B_2,\mathcal{O}_{B_2}(H-E_2))\\
&=&H^i(B_1,\mathcal{I}_{\widetilde{\Sigma}}(H))\\
&=&H^i(\PP^n,\mathcal{I}_{\Sigma}(1)),
\end{eqnarray*}
where the last equality is a consequence of $\Sigma$ having rational singularities.
From the sequence on $B_2$
$$0\rightarrow\mathcal{O}_{B_2}(H-E_1-E_2)\rightarrow\mathcal{O}_{B_2}(H-E_1)\rightarrow\mathcal{O}_{E_2}(H-E_1)\rightarrow0$$
Once again considering $\widetilde{\Sigma}$ as a $\PP^1$-bundle over $\PP^2$, we see that $\mathcal{O}_{\widetilde{\Sigma}}(H-E_1)$ is $\mathcal{O}_{\PP^1}(-1)$ along the fibers, thus $H^i(\widetilde{\Sigma},\mathcal{O}_{\widetilde{\Sigma}}(H-E_1))=0$ for $i\geq0$.  Putting this together, we have
\begin{eqnarray*}
H^i(\PP^n,\mathcal{I}_{\Sigma}(1))&=&H^i(B_2,\mathcal{O}_{B_2}(H-E_1-E_2))\\
&=&H^i(B_1,\mathcal{O}_{B_1}(H-E_1))\\
&=&H^i(\PP^n,\mathcal{I}_{X}(1))\\
&=&0.
\end{eqnarray*}

For $i=1$, in a similar fashion it is enough to show $H^1(B_2,\mathcal{O}_{B_2}(2H-E_1-E_2))=0$.  From the sequence
$$0\rightarrow\mathcal{O}_{B_2}(2H-E_1-E_2)\rightarrow\mathcal{O}_{B_2}(2H-E_1)\rightarrow\mathcal{O}_{E_2}(2H-E_1)\rightarrow0$$
and the fact that $H^1(B_2,\mathcal{O}_{B_2}(2H-E_1))=H^1(\PP^n,\mathcal{I}_X(2))=0$, it suffices to show $H^0(B_2,\mathcal{O}_{B_2}(2H-E_1))\rightarrow H^0(E_2,\mathcal{O}_{E_2}(2H-E_1))$ is surjective.  However, we know that $H^0(B_2,\mathcal{O}_{B_2}(2H-E_1-E_2))=H^0(B_1,\mathcal{I}_{\widetilde{\Sigma}}(2))=0$, which show that the map is injective.  Thus, if we shows the two spaces have the same dimension we are done.   We have the well-known identification $H^0(B_2,\mathcal{O}_{B_2}(2H-E_1))=H^0(\PP^n,\mathcal{I}_X(2))$.  Further, letting $\varphi:B_1\rightarrow\PP^s$ be the map given by quadrics vanishing on $X$ (i.e. the morphism induced by the linear system $|2H-E_1|$), then we have the restriction $\overline{\varphi}:\widetilde{\Sigma}\rightarrow\PP^2$.  

Note that $\widetilde{\Sigma} \subset \PP^n \times \PP^s$ is a $\PP^1$-bundle over $\PP^2$.  It is, further, a nice exercise to show that the double cover $\PP^1\times\PP^1=\widetilde{\Sigma}\cap E_1\rightarrow\PP^2$ is, in situ, the natural double cover re-embedded by the Veronese $v_{n-2}$.  Therefore, $\mathcal{O}_{\PP^s}(1)|_{\PP^2}=\mathcal{O}_{\PP^2}(n-2)$, and this implies that $\overline{\varphi}^*\mathcal{O}_{\PP^2}(n-2)=\mathcal{O}_{\widetilde{\Sigma}}(2H-E_1)$.  Hence, $H^0(\widetilde{\Sigma},\mathcal{O}_{\widetilde{\Sigma}}(2H-E_1))=H^0(\PP^2,\mathcal{O}_{\PP^2}(n-2))$.  A quick computation gives $h^0(\PP^n,\mathcal{I}_X(2))=h^0(\PP^2,\mathcal{O}_{\PP^2}(n-2))=\binom{n}{2}$.
\end{proof}

In fact, we see from the proof that  $H^1(\PP^n,\mathcal{I}_{\Sigma}(1))=0$.  As 3-regularity implies that $H^1(\PP^n,\mathcal{I}_{\Sigma}(k))=0$ for all $k \geq 2$ we also have:
\begin{cor}
$\Sigma$ is projectively normal.
\qed
\end{cor}

\appendix
\section{Code for computing prolongations}
The code in this section was written to produce examples for joint work of the first author and Seth Sullivant.  Mike Stillman wrote the intersection.m2 package and helped rewrite the code for prolongations.

The code in intersection.m2 was written to allow us to specify a degree $d$ and intersect the degree $d$ piece of a list of ideals.
\begin{verbatim}
-- file name: intersection.m2
-- intersect needs a degree limit

intersection = method(Dispatch => Thing, 
TypicalValue=>Ideal, Options => {DegreeLimit=>{}})

intersection Sequence := 
intersection List := o -> L -> (
     if not all(L, x -> instance(x,Ideal))
     then error "expected a list of ideals";
     B := directSum apply(L, generators);
     A := map(target B, 1, (i,j) -> 1);
     ideal syz(A|B, SyzygyRows => 1, DegreeLimit=>o.DegreeLimit)
     )

--INPUT: L = list of ideals
--       d = integer
--OUTPUT: generators for the intersection of the ideals in L in degrees <=d
inter = (d, L) -> (
     ans := L_0;
     for i from 1 to #L-1 do time ans = intersection(ans,L_i,DegreeLimit=>d);
     ans)

--INPUT:  L = list of ideals, each generated in degree d
--        d = integer
--OUTPUT:  the ideal generated by the intersection of the degree d parts
 of each ideal

intersectSpaces = (d,L) -> (
     S := ring L_0;
     k := coefficientRing S;
     monoms := basis(d,S);
     time L = apply(L, I -> (
	       (mn,cf) := coefficients(gens I, Monomials=>monoms);
	       image substitute(cf,k)));
     time M := intersect L;
     ideal(monoms * (gens M))
     )
\end{verbatim}

The code for computing prolongations follows below.  If $X \subset \PP^n$ has ideal $I$ generated by quadrics, then the prolongation $P$ of a basis for the quadrics in $I$ contains all forms of degree 3 in $I(\Sigma_1(X)).$  Of course, $P$ is not guaranteed to generate the ideal of the secant variety.  In practice, it is often easy to verify that the ideal of the secant variety is generated by $P$ by computing the projective dimension, degree, and dimension of $S/\langle P \rangle.$  If we see that the degree and dimension are as expected and $S/\langle P \rangle$ is Cohen-Macaulay, then we conclude that $I(\Sigma_1(X)) = \langle P \rangle.$

\begin{verbatim}
load "intersection.m2"

--INPUT: f = polynomial in S
--   	   i = index of i-th variable in S
--OUTPUT: an antiderivative of f integrated with respect to the i-th 
variable in S.

integration = (i,f) ->(
     S:=ring f;
    sum apply (terms f, g -> (1/((flatten exponents g)_i+1)) * g*S_i)
     )

--INPUT: A = a list of forms of degree d.
--OUTPUT: a basis for the prolongation of A.
prolong = (A) ->(
     S:= ring A_0;
     n:= numgens S;
     d:= (degree(A_0))_0;
     L :={};
     
     for i from 0 to (n-1) do(
	 M := ideal drop (flatten entries vars S, {i,i});
	 intA:=ideal apply(A, f -> integration(i, f))+M^(d+1);
	 L = append(L, intA);
	  );
     << "d = " << d << endl;
     LL = L;
     GenList := flatten entries mingens intersectSpaces(d+1,L) 
)
\end{verbatim}

\end{document}